\documentstyle[12pt,amssymb]{article}
\begin{document}
\newtheorem{theorem}{Theorem}
\newtheorem{lemma}{Lemma}
\newtheorem{definition}[theorem]{Definition}
\newtheorem{corollary}[theorem]{Corollary}
\newtheorem{proposition}[theorem]{Proposition}
\newcommand{\bproof}{\noindent{\bf Proof}}
\newcommand{\kproof}{\noindent{\bf Proof:}}
\newcommand{\eproof}{\rule{0pt}{1ex}\hfill $\Box$\par}
\newtheorem{remark}{Remark}
\setcounter{section}{0}

\title{Vector-valued $L_p$-convergence of orthogonal series and Lagrange 
interpolation.}
\author{H. K\"onig (Kiel)\thanks{Supported in part by the DFG, Ko
962/3-1.}\and N. J. Nielsen (Odense)\thanks{Supported in part by the DFG, Ki
51/98-1 and the Danish Natural Science Council, grant 11-8622.}}
\date{ }
\maketitle
\vspace{1cm}

\begin{abstract}
We give necessary and sufficient conditions for interpolation
inequalities of the type considered by Marcinkiewicz and Zygmund to be true 
in the case of Banach space-valued polynomials and Jacobi weights and nodes. 
We also study the vector-valued expansion problem of $L_p$-functions in terms 
of Jacobi polynomials and consider the question of unconditional
convergence. The notion of type $p$ with respect to orthonormal systems
leads to some characterizations of Hilbert spaces. It is also shown that
various vector-valued Jacobi means are equivalent.
\end{abstract}

\section{Introduction and results}

Let $X$ be a Banach space, $1 \le p < \infty$ and $L_p(\Bbb R;X)$
denote the space of (classes of) p-th power integrable functions with
norm $\parallel f \parallel := (\int_{\Bbb R} \parallel f (t)
\parallel^pdt)^{1/p}$. A Banach space is a {\em UMD-space} provided
that the Hilbert transform on $\Bbb R$,

\begin{equation}
Hf(t):= \quad \mbox{p.v.} \int_{\Bbb R} \ \frac{f(s)}{t-s}ds, \  f \in
L_p(\Bbb R;X),
\end{equation}

\noindent defines a bounded operator $H:\ L_p(\Bbb R;X) \longrightarrow L_p(\Bbb
R;X)$ for some $1<p<\infty$. It is well-known that this holds for {\em
some} $1<p<\infty$ if and only if it holds for {\em all} $1<p<\infty$, see
e.g. Schwarz \cite{kn:sch}. All $L_q(\mu)$-spaces with $1<q<\infty$
or all reflexive Orlicz spaces are UMD-spaces, cf. Fernandez and Garcia
\cite{kn:fg}.

Let $I=(-1,1)$, $\alpha, \beta>-1$ and $w_{\alpha
\beta}(t):=(1-t)^\alpha(1+t)^\beta$ for $t \in I$. Let 

\[
L_p(I,w_{\alpha\beta};X):=\{ f:I \longrightarrow X | \parallel f
\parallel _p:= \parallel f \parallel _{p;\alpha , \beta}:=(\int_I
\parallel f(t) \parallel ^p w_{\alpha \beta} (t) dt)^{1/p} < \infty \}.\]

The scalar product in $L_2(I,w_{\alpha \beta}):=L_2 (I,w_{\alpha
\beta};\Bbb R)$ will be denoted by $<\cdot,\cdot>$ or $<\cdot,\cdot
>_{\alpha \beta}$. For $\alpha = \beta$ we just write $w_\alpha$ and $<\cdot
,\cdot>_\alpha$. By $\Pi_n(X)$ we denote the space of polynomials of 
degree $\le n$ with coefficients in $X$. Let $\Pi_n := \Pi_n (\Bbb R)$. The
$L_2(I,w_{\alpha \beta})$-normalized Jacobi polynomials with respect to
$(I, w_{\alpha \beta})$ will be denoted by $p^{(\alpha ,\beta)}_n $, $n
\in \Bbb N_0$. Hence $p_n^{(\alpha ,\beta)} \in \Pi_n$ and

\begin{equation}
<p_n^{(\alpha , \beta)}, p_m^{(\alpha , \beta)} >_{\alpha , \beta} =
\int_I p_n^{(\alpha , \beta)}(t)p_m^{(\alpha , \beta)}(t)
w_{\alpha \beta} (t) dt \ = \ \delta_{nm}.
\end{equation}

This normalization is more convenient for us than the standard one of
Szeg\"{o} \cite{kn:sz}. For $\alpha = \beta = -\frac{1}{2}$
 $(\frac{1}{2})$ one gets the Tchebychev polynomials of the first (second)
kind, for $\alpha = \beta = 0$ the Legendre polynomials. Let $t_1 >
\cdots > t_{n + 1}$ denote the zeros of $p^{(\alpha, \beta)}_{n+1}$, all
of which are in $I$, and $\lambda_1, \cdots , \lambda _{n + 1} > 0$ the
Gaussian quadrature weights. Thus for any real polynomial $q$ of degree
$\le 2n + 1$, one has

\begin{equation}
\int_I q(t) w_{\alpha \beta}(t) dt \quad = \quad \sum_{j=1}^{n+1}
 \lambda_j q (t_j).
\end{equation}

Clearly, $\lambda_j$ and $t_j$ depend on $n, j, \alpha$ and $\beta$ but
not on $q$. One has for $\alpha, \beta > -1$

\begin{equation}
\lambda_j = (2n + \alpha + \beta + 3)((1-t_j^2)
p_{n+1}^{(\alpha, \beta)'}(t_j)^2)^{-1} \quad \sim 
\quad \left\{\begin{array}{rl}
{j^{2\alpha + 1}}/{n^{2\alpha +2}} & \mbox{$j \le \frac{n}{2}$} \\
{(n+2-j)^{2\beta+1}}/{n^{2\beta+2}} & \mbox{$j > \frac{n}{2}$}
\end{array}\right\}
\end{equation}

\begin{eqnarray}
1 - t^2_j \sim (j/n)^2, and & p_{n+1}^{(\alpha ,
\beta)'} (t_j) \sim {n^{5/2+\alpha}}/{j^{3/2 + \alpha}} &
\mbox{for all $j \le \frac{n}{2}$}.\\
  & p^{(\alpha , \beta)}_n (-x) = 
(-1)^n p_n^{(\alpha , \beta)}(x). & \mbox{} \nonumber
\end{eqnarray}

See Szeg\"{o} \cite[3.4, 4.1, 4.3, 8.9, 15.3]{kn:sz}, taking into
account the different normalization there. Here $\lambda_j \sim f_j$
means that there are constants $c_1, c_2 > 0$ independent of $j$ and $n$
such that $c_1 f_j \ \le \ \lambda_j \ \le \ c_2 f_j$ for all $n$ and
$j$ concerned. For $\alpha = \beta = -\frac{1}{2}$, $\lambda_j =
\pi/(n+1)$.

Marcinkiewicz and Zygmund \cite[ch. X]{kn:z1} proved interpolation
inequalities for trigonometric polynomials of degree $\le n$ which for
even trigonometric polynomials $g$, after a transformation
$g(x)=q(\mbox{cost}), \ x=\mbox{cost}, \ q \in \Pi_n$, can be restated
as

\[
\frac{1}{3} (\sum^{n+1}_{j=1} |q(t_j)|^p /(n+1))^{\frac{1}{p}} \ \le
 \ (\frac{1}{2} \int_I |q(t)|^p (1-t^2)^{-\frac{1}{2}}dt)^{\frac{1}{p}} \ \le
 \ c_p(\sum^{n+1}_{j=1} |q(t_j)|^p/(n+1))^{1/p}.
\]

Here $(t_j)$ are the zeros of the Jacobi polynomial $p_{n+1}^{\alpha ,
\beta}$ in the Tchebychev case $\alpha = \beta = -1/2$, and $c_p$
depends on $1<p<\infty$ only. The left inequality holds for $p=1,
\infty$ as well whereas the right one fails, in general. For $p=2$,
(3) gives more precise information since $\lambda_j = \pi /(n+1)$.
The Marcinkiewicz-Zygmund inequalities extend to the Jacobi case of
general $\alpha, \beta > -1$ and to the vector-valued setting in the
following sense:

\begin{theorem}
Let $X$ be a Banach space, $1 \le p \le \infty$, $\alpha, \beta > -1$,
 $(t_j)$ the zeros of $p_{n+1}^{\alpha , \beta}$ and $(\lambda_j)$ the
corresponding quadrature weights.

\begin{itemize}
\item [a)] There is $c > 0$ such that for all $1 \le p \le \infty$, $n
\in \Bbb N$ and $q \in \Pi_{2n}(X)$

\begin{equation}
c^{-1}(\sum_{j=1}^{n+1} \lambda_j \parallel q(t_j)\parallel^p)^{1/p} \le
(\int^1_{-1} \parallel q(t) \parallel ^p w_{\alpha \beta}(t) dt)^{1/p}.
\end{equation}

\item [b)] Let 
\begin{eqnarray*}
\mu (\alpha ,\beta): & = & \max (1,4 (\alpha +1)/(2\alpha
+5), 4 (\beta + 1)/(2\beta + 5))\\
 m(\alpha, \beta): & = & \max
(1,4(\alpha +1)/(2\alpha +3), 4 (\beta +1)/(2\beta +3))
\end{eqnarray*}
and
$M(\alpha, \beta):= m(\alpha, \beta)'$, i.e. $m(\alpha, \beta)^{-1} +
M(\alpha, \beta)^{-1} =1$. Then the following are equivalent. 

\begin{enumerate}
\item [(1)] There is $c_p > 0$ such that for all $n \in \Bbb N$ and $q
\in \Pi_n(X)$

\begin{equation}
(\int_{-1}^1 \parallel q(t) \parallel^p w_{\alpha \beta}(t) dt)^{1/p} 
\le c_p (\sum_{j=1}^{n+1} \lambda_j \parallel q (t_j) \parallel^p)^{1/p}.
\end{equation}

\item [(2)] $X$ is a UMD-space and $p$ satisfies $\mu (\alpha, \beta) < p <
M(\alpha, \beta)$.
\end{enumerate}
\end{itemize}
\end{theorem}

Part (a) is proved just as the scalar result which goes back to Askey
\cite{kn:as}, Nevai \cite{kn:n} and Zygmund \cite{kn:z2}. The converse 
inequality
(7) was shown in the scalar case (for $\alpha = \beta$) by Askey
\cite{kn:as} under the more restrictive assumption $m(\alpha, \beta) < p <
M(\alpha, \beta)$ using (a) and duality; the duality method, however,
fails if $\mu(\alpha, \beta) < p \le m(\alpha, \beta)$. The question
whether (7) in the vector-valued case requires $X$ to be a UMD-space
was raised by Pietsch in the case of trigonometric polynomials
(corresponding to $\alpha = \beta = -1/2$) and solved by him in this
case by a different method \cite{kn:pie}.

In terms of Banach spaces, Theorem 1 states that the spaces $ \Pi_n(X)_p$ 
as subspaces of $L_p(I,w_{\alpha \beta} ;X)$ 
are uniformly isomorphic to $l^{n + 1}_p(X)$-spaces, by evaluating 
the polynomials $q$ at the zeros $(t_j)$, provided that (b), (2) holds;
i.e. the Banach-Mazur distances $d(\Pi_n(X)_p,l_p^{n+1}(X))$ are uniformly
bounded.

For $f \in L_p (I, w_{\alpha \beta} ; X)$, let $Q_nf:= \ \sum^n_{j=0} <
f,p_j^{\alpha,\beta}> p_j^{(\alpha, \beta)} \in \Pi_n(X)_p$ denote the
orthogonal projection of $f$ onto the space of polynomials of degree $\le n$. 
The
following vector-valued expansion theorem for Jacobi polynomials
generalizes the classical scalar result of Pollard \cite{kn:pol} and
Muckenhaupt \cite{kn:mu}.

\begin{theorem}
Let $X$ be a Banach space, $1 \le p \le \infty$, $\alpha, \beta > -1$ and
$m(\alpha, \beta)$ and $M(\alpha, \beta)$ as before. Then the following
are equivalent:

\begin{itemize}
\item [(1)] For all $f \in L_p (I,w_{\alpha \beta}; X) \quad Q_nf$
converges to $f$ in the $L_p$-norm.
\item [(2)] $X$ is a UMD-space and $m(\alpha, \beta) < p < M(\alpha,
\beta)$.
\end{itemize}
\end{theorem}

The necessity of the UMD-condition on $X$ will be proved using Theorem
1; the interval for $p$ is ``symmetric'' with respect to $p=2$ and smaller than
the one exhibited in Theorem 1, (b). Analogues of Theorems 1 and 2 in the
case of the Hermite polynomials are proved in \cite{kn:ko}. Using the
results of Gilbert \cite{kn:gi}, we also prove that various vector-valued
Jacobi means are equivalent:

\begin{proposition}
Let $\alpha, \beta > -1$, $1< p < \infty$, $\gamma \in \Bbb R$ with
$|\frac{\gamma}{2} + \frac{1}{p} - \frac{1}{2} | < \frac{1}{4}$. Let $X$
be a UMD-space. Then there is $M=M(\alpha, \beta, \gamma , p) \ge 1$
such that for all $n \in \Bbb N$ and all $x_1, \dots ,x_n \in X$

\begin{eqnarray}
& & (\int^1_{-1} \parallel \sum^n_{j=0} p_j^{(\alpha, \alpha)} (t)x_j
\parallel^p (1-t^2)^{(\alpha + \gamma)p/2} dt)^{1/p} \nonumber \\ 
& \stackrel{M}{\sim} & 
 (\int^1_{-1} \parallel \sum^n_{j=0} p_j^{(\beta, \beta)}(t)x_j
\parallel^p (1-t^2)^{(\beta + \gamma)p/2} dt)^{1/p}. 
\end{eqnarray}
\end{proposition}

Here $\stackrel{M}{\sim}$ means that the quotient of the two expressions is
between $1/M$ and $M$. Instead of $(\alpha, \alpha)$ and $(\beta,
\beta)$, one could consider $(\alpha_1, \alpha_2)$ and $(\beta_1,
\beta_2)$ as Jacobi-indices, provided the weight functions are changed
accordingly. The convergence of the Jacobi series in Theorem 2 is not
unconditional unless $p=2$ and $X$ is a Hilbert space, as will follow
from the following general result. Recall that a series $\sum_{n \in
\Bbb N}y_n$ in a Banach space $Y$ converges {\em unconditionally} if
$\sum_{n \in \Bbb N} \varepsilon_ny_n$ converges in $Y$ for all choices
of signs $\varepsilon_n = \pm 1$.

\begin{proposition}
Let $(\Omega, \mu)$ be a measure space and $(p_n)$ be a complete
orthonormal system in $L_2(\Omega, \mu)$, assumed to be infinite 
dimensional. Let $X$ be a Banach space and $1 \le p < \infty$. Assume
that for all $f \in L_p(\Omega, \mu; X)$, the series $\sum_n < f,p_n
>p_n$ converges unconditionally in $L_p (\Omega, \mu; X)$. Then:

\begin{itemize}
\item [(i)] If $\parallel p_j \parallel_2 \sim \parallel p_j
\parallel_{\max (p,p')}$ and $(\Omega, \mu)$ is a finite measure space,
one has $p=2$.
\item [(ii)] If $\mbox{sup}_j |p_j| \in L_2(\Omega, \mu), \ X$ is
isomorphic to a Hilbert space.
\end{itemize}
\end{proposition}

Statement (ii) was also shown by Defant and Junge \cite{kn:dj}. Both
conditions (i) and (ii) are satisfied in the Jacobi case provided that
the condition $m(\alpha,\beta) < p < M(x,\beta)$ holds (necessary for
convergence). Without an assumption like $\mbox{sup}_j |p_j|\in L_2
(\Omega, \mu)$, $X$ is not isomorphic to a Hilbert space in general, as
the Haar system shows. However, one has:

\begin{proposition}
Let $1<p<\infty$ and $(p_n)_{n \in \Bbb N}$ be an unconditional basis of
$L_p(0,1)$. Let $X$ be a Banach space such that for any $f \in
L_p(0,1;X)$, the series $\sum_{n \in \Bbb N} <f,p_n>p_n$ converges
unconditionally in $L_p(0,1;X)$. Then $X$ is a UMD-space.
\end{proposition}

The proof shows that the Haar basis is unconditional in $L_p(0,1;X)$
which by Maurey \cite{kn:ma}, Burkholder \cite{kn:bu} and Bourgain
\cite{kn:bo} is equivalent to $X$ being a UMD-space. It was shown by
Aldous \cite{kn:al} that $X$ is a UMD-space if $L_p(X)$ has an
unconditional basis.

Let $(\Omega, \mu)$ be a measure space and $(p_n)_{n \in \Bbb N}$ be a
complete orthonormal system in $L_2(\Omega, \mu)$. We say that a Banach
space $X$ has ${\bf (p_n) -type}$ $2$ provided there is $c>0$ such that
for all $m \in \Bbb N$ and all $x_1, \dots ,x_m \in X$

\[
(\int_{\Omega} \parallel \sum^m_{j=1} p_j(t) x_j \parallel^2
d\mu(t))^{1/2} \le c (\sum^m_{j=1} \parallel x_j \parallel^2)^{1/2}.
\]
 
$X$ has ${\bf (p_n) -cotype}$ $2$ if the reverse inequality holds. In
\cite{kn:pi2}, Pisier showed for the Haar system $(h_n)$, that
$(h_n)$-type 2 of $X$ is equivalent to $X$ being 2-smooth, e.g. has an 
equivalent uniformly convex norm with modulus of convexity of power type 2. In
Pisier and Xu \cite{kn:px} the related notion of $H$-type $p$ $(\le 2)$ is
considered for all orthonomal systems $(p_n)$. Kwapie\'n \cite{kn:kw}
studied this notion for the trigonometric system $(e_n)$ in
$L_2(0,2\pi)$ showing that $(e_n)$-type 2 (also called Fourier-type 2,
$e_n(t) = exp(int)$) of $X$ implies that $X$ is isomorphic to a
Hilbert space. This result generalizes to the case of Jacobi
polynomials.

\begin{proposition}
Let $X$ be a Banach space which is Jacobi $(p_n^{(\alpha ,\beta)})$-type
2 for some $\alpha ,\beta > -1$. Then $X$ is isomorphic to a Hilbert
space.
\end{proposition}

The proof uses the interpolation inequalities of Theorem 1. In general,
$(p_n)$-type 2 implies type 2 in the usual sense \cite{kn:mp} , i.e.
with respect to the Rademacher system $(r_n)$, $r_n(t)= \ \mbox{sgn sin}
2^n \pi t$.

\begin{proposition}
Let $X$ be a Banach space which is of $(p_n)$-type
2 for some complete orthonormal system $(p_n)$ in $L_2(0,1)$. Then $X$
is of Haar type 2, hence 2-smooth and of type 2.
\end{proposition}

There is a partial converse to this result.

\begin{proposition} Let $X$ be a Banach space and $(p_n)_{n \in \Bbb N}
\subset L_2(0,1)$ be a complete orthonormal system such that for any $f
\in L_2(0,1;X)$ the series $\sum_{n \in \Bbb N} < f, p_n > p_n$
converges unconditionally in $L_2(0,1;X)$. Then, if $X$ has type 2, it
also has $(p_n)$-type 2.
\end{proposition}

It follows from Proposition 5 that the unconditionality assumption in
Proposition 8 implies that the space $X$ in question has UMD. On the other 
hand, if $X$ has UMD the unconditionality assumption in Proposition 8 is
satisfied for the Haar system, and thus by Pisier's result mentioned above,
type 2 and UMD of $X$ implies that $X$ is 2-smooth. This in turn implies 
type 2 but does not imply the UMD-property, since by Bourgain \cite{kn:bo} 
there exists a Banach lattice satisfying an upper-p and lower-q estimate and 
failing the UMD-property; choosing $2<p<q<\infty$ there, such a lattice is 
2-smooth, cf. \cite{kn:lt}.

\section{The interpolation inequalities}
For the proof of Theorem 1, we need a well-known fact about continuity
in $L_p$, cf. Pollard \cite{kn:pol} or Benedek, Murphy and Panzone
\cite{kn:bmp}. In the scalar case, it is a special case of the theory 
of weighted singular integral operators with weights in the Muckenhaupt
class $A_p$, cf. Garcia-Cuerva and de Francia \cite[chap. IV]{kn:gr}.

\begin{lemma} Let $X$ be a Banach space, $1 \le p \le \infty, \ b \in
\Bbb R$ and $k:\Bbb R^2 \longrightarrow \Bbb R$ be defined by $k(u,v):=
\ | |u/v|^b -1| / |u-v|$. Then the integral operator $T_k$ given by
$T_kf(u):= \int_{\Bbb R} k(u,v)f(v)dv$ defines a bounded operator
$T_k:L_p(\Bbb R;X) \longrightarrow L_p (\Bbb R;X)$ provided that $-1/p <
b < 1-1/p$ (actually if and only if).
\end{lemma}

{\kproof}
We sketch the simple proof. Let $r(u,v):= |u/v|^{1/pp'}$. It suffices
to show that 

\begin{equation}
\sup_u \ \int_{\Bbb R} \ k(u,v) \ r(u,v)^{p'} dv \le M, \ \sup_v \
\int_{\Bbb R} \ k(u,v) r(u,v)^{-p} du \ \le \ M.
\end{equation}

An application of H\"older's inequality then shows that $T_k$ is
continous as a map \mbox{$T_k:L_p(\Bbb R; X) \longrightarrow L_p(\Bbb R; X)$} 
with norm $\le M$. To check the first inequality in (9), substitute $v/u
=t$ to find

\[
\sup_{u \ne 0} \int_{\Bbb R} k(u,v) \ r(u,v)^{p'}dv \ = \ \int_{\Bbb R}
| t^{-b} -1| \ |t|^{-1/p}/|t-1| dt.
\]

This is finite since integrability at 0 is assured by $b < 1-1/p$,
and integrability at $\pm \infty$ by $b > -1/p$. Note that for $t
\longrightarrow 1$, there is no singularity, the integrand tends to
$|b|$. The second condition in (9) is checked similarly.
{\eproof}

By Szeg\"{o} \cite[7.32 , 4.3]{kn:sz}, for any pair of indices $\alpha,
 \beta > -1$, there is $c=c_{\alpha, \beta}$ such that for all $n \in
\Bbb N$ and $t \in [-1,1]$, the $L_2$-normalized Jacobi polynomials
$p_n^{(\alpha ,\beta)}$ satisfy the estimate

\begin{equation}
|p_n^{(\alpha ,\beta)} (t)| \le \ c \ (1-t+n^{-2})^{-(\alpha /2 + 1/4)}
(1+t+n^{-2})^{-(\beta /2 + 1/4)}
\end{equation}

{\bproof} {\bf of theorem 1:} We start with 

(b), (2)$\Rightarrow$ (1).

\noindent Assume that $X$ is a UMD space and that $p$
satisfies $\mu(\alpha ,\beta) < p < M(\alpha ,\beta )$. Let $q \in
\Pi_n(X)$ and put $y_j:= q(t_j)/p'_{n+1}(t_j)$. The Lagrange functions
$\ell_j \in \Pi_n$,

\[
\ell_j(t):= \ p_{n+1}^{(\alpha ,\beta)} (t) \ / \ ({p_{n+1}^{(\alpha ,
\beta)}}'(t_j)(t-t_j))
\]
satisfy $\ell_j(t_i)=\delta_{ji}$ for $i,j=1,\cdots, n+1$ and thus $q$
coincides with its interpolating polynomial $q = \sum_{j=1}^{n+1}
q(t_j)\ell_j$. We have to estimate

\[
L:= ( \int^1_{-1}\ \parallel q(t) \parallel^p w_{\alpha \beta}(t)
dt)^{1/p} \ = \ (\int^1_{-1} \ \parallel \sum^{n+1}_{j=1} y_j \
\frac{p_{n+1}^{(\alpha , \beta)} (t)}{t - t_j} \ \parallel^p w_{\alpha
\beta} (t) dt)^{1/p}
\]
from above. Let $I_j = (t_j, t_{j-1}), \ |I_j| = (t_{j-1}- t_j)$ and
$\chi_j$ be the characteristic function of $I_j$, for $j=1, \cdots , n+1$,
with $t_o:=1$. The proof relies on the fact that $1/(t-t_j)$ is
sufficiently close to the Hilbert transform of $-\chi_j/|I_j|$ at $t$ which
is

\[ 
H(- \frac{\chi_j}{|I_j|})(t) = \frac{1}{|I_j|} \log| \frac{t -
t_{j-1}}{t-t_j}| = \frac{1}{|I_j|} \log|1 - \frac{|I_j|}{t-t_j}|.
\]

Let $J_n = [a_n ,b_n]$ where
\[
a_n=\left\{\begin{array}{cl}
-1 & \mbox{if $\beta > -1/2$} \\
-1 + dn^{-2} & \mbox{if $\beta \le -1/2$},
\end{array} \right.
\]

\[
b_n=\left\{\begin{array}{cl}
1 & \mbox{if $\alpha > -1/2$} \\
1 - dn^{-2} & \mbox{if $\alpha \le -1/2$;}
\end{array} \right.
\]

\noindent and $d$ is chosen such that $\min (1-t_1, 1 + t_{n+1}) \ge 2d n^{-2}$. 
By
\cite{kn:sz} this is possible. 

It follows from (10) that for $n \in \Bbb N$ and $t \in J_n$

\begin{equation}
|p_n^{(\alpha, \beta)}(t)| \le c(1-t)^{-(\alpha/2 + 1/4)}
(1+t)^{-(\beta/2 + 1/4)}.
\end{equation}

In the following, constants $c_1, c_2, \cdots $ may depend on $\alpha,
\beta$ and $p$, but not on $n, j$ and $t$. We claim that for $n \in \Bbb
N$ and $t \in J_n$

\begin{equation}
|p_{n+1}^{(\alpha,\beta)} (t)|| \frac{1}{t-t_j} +
H(\frac{\chi_j}{|I_j|})(t) | \le f_j(t),
\end{equation}

\noindent where

\[
f_j(t) := c_1 \min (\frac{1}{|I_j|}, \frac{|I_j|}{(t-t_j)^2})
(1-t)^{-(\alpha/2 + 1/4)}(1+t)^{-(\beta/2 + 1/4)}.
\]

If $t$ is such that $|t-t_j| > 2|I_j|$, (12) follows from (11) and $|x -
\log(1+x)|\le x^2$ for $|x| \le 1/2$, i.e. $|\frac{1}{t-t_j} +
H(\frac{\chi_j}{|I_j|})(t)| \le \frac{|I_j|}{(t-t_j)^2}$. 
For $|t-t_j| \le 2 |I_j|$, one uses that $p_{n+1}^{(\alpha, \beta)}$ has
a zero in $t_j$. By the mean-value theorem there is a $\theta$ between
$t$ and $t_j$ such that 
\[
|\frac{p_{n+1}^{(\alpha,\beta)}(t)}{t-t_j} |= |\frac{p_{n+1}^{(\alpha,
\beta)} (t) - p_{n+1}^{(\alpha ,\beta)}(t_j)}{t-t_j} | = |
{p_{n+1}^{(\alpha ,\beta)}}'(\theta) | \le f_j (t),
\]
using that by Szeg\"o \cite[8.9]{kn:sz} and (5), e.g. for $j \le n/2$,
\[
|{p_{n+1}^{(\alpha ,\beta)}}'(\theta) | \le c_3 n^{\alpha +
5/2}/j^{\alpha + 3/2} \sim f_j (t_j) \sim f_j (t).
\]

We note that (only) for $j=1$, the logarithmic singularity  of
$H(x_j/|I_j|)$ at $t=1$ is {\em not} compensated by a zero of
$p_{n+1}^{(\alpha ,\beta)}$ $(t_0=1)$, but (12) is true in this case and $
\alpha > -1/2$ too, since by (10) for $(1-t) \le n^{-2}$

\begin{eqnarray*}
|p_n(t) H(\frac{\chi_j}{|I_j|}) (t) | & \le & c_4 n^{5/2 +
\alpha}|\log(n^2(1-t))| \\
 & \le & c_5 n^2 (1-t)^{-(\alpha/2 + 1/4)}
\end{eqnarray*}

\noindent using $|\log v| \le c_3 v^{-\varepsilon}$ for $\varepsilon =
\frac{\alpha}{2} + \frac{1}{4} > 0$ and $0 < v \le 1$. Hence (12) holds.
Applying this we find

\[
L_n := (\int_{J_n} \parallel \sum^{n+1}_{j=1} y_j
\frac{p_{n+1}^{(\alpha ,\beta)} (t)}{t-t_j} \parallel^p w_{\alpha \beta}
(t) dt )^{1/p} \le M_1 + M_2,
\]
where

\begin{eqnarray*}
M_1 & := & c_6 (\int_{J_n} \parallel H(\sum_{j=1}^{n+1} y_j \
\chi_j/|I_j|)(t) \parallel^p (1-t)^{\gamma p}(1+t)^{\delta p} dt)^{1/p},\\
M_2 & := & (\int_{J_n} (\sum_{j=1}^{n+1} \parallel y_j \parallel f_j
(t))^p w_{\alpha \beta} (t) dt)^{1/p},
\end{eqnarray*}
and $\gamma := \alpha(1/p - 1/2) - 1/4$, $\delta := \beta(1/p - 1/2)
-1/4$. The restrictions on $p$ are equivalent to $-1/p < \gamma$, $\delta
< 1 - 1/p$ and $1 < p < \infty$. In particular, $|\gamma | < 1$, $|\delta
| < 1$. We estimate the ``main'' term $M_1$ and the ``error'' term $M_2$
separately.

We claim that the kernel \mbox{$K(t,s) := 1/(t-s)
((1-t)/(1-s))^{\gamma} ((1+t)/(1+s))^{\delta}$}, $ t, s \in [-1,1]$
defines a bounded integral operator $T_k : L_p (-1,1; X) \longrightarrow
L_p (-1,1;X)$. Indeed, by lemma 1, the kernel \mbox{$|((1-t)/(1-s))^{\gamma} -1 
|
/ |t-s|$} defines a bounded operator $L_p (\Bbb R; X) \longrightarrow L_p
(\Bbb R;X)$, replacing $t$ and $s$ by $(1-t)$ and $(1-s)$. Since $X$ is
a UMD-space, so does $1/(t-s)$ and hence also \mbox{$1/(t-s)
((1-t)/(1-s))^{\gamma}$}. Since \mbox{$(\frac{1+t}{1+s})^{\delta}$} is
bounded from above and below by positive constants for $s, t \in [0,1]$, the
kernel $K$ defines a bounded operator $T_k : L_p (0,1; X)
\longrightarrow L_p (0,1; X)$. The same holds on the interval $[-1,0]$.
The kernel is less singular for $t$ and $s$ of different sign: if $s \in
[-1,0]$, $t \in [0,1]$, the substitution $s \rightarrow -s$ yields a
kernel of the type \mbox{$(1/(t+s)) w_1 (t) w_2 (s)$} on $[0,1]^2$, where
$w_1$ and $w_2$ are integrable over $[0,1]$ and bounded near $0$. By
Hardy, Littlewood and Polya \cite{kn:hlp}, the kernel $1/(t+s)$ defines
an operator $L_p (0, \infty ) \longrightarrow L_p (0, \infty )$ of norm
$\pi / sin(\pi /p)$, for any $1 < p < \infty$. Since $1/(t+s)$ is
positive, this also holds for $X$-valued functions and hence $T_k :
L_p(-1,0; X) \longrightarrow L_p (0,1; X)$ is bounded as well. The case
$s \in [0,1]$, $t \in [-1,0]$ is treated similarly. Together these facts
prove the claim. Hence there is a $c_7$ such that for all $f \in L_p (-1,
1; w_{\gamma \delta}; X)$
\[
(\int_{-1}^1 \parallel \int_{-1}^1 \frac{f(s)}{t-s} ds \parallel^p
(1-t)^{\gamma p} (1+t)^{\delta p} dt)^{1/p} \le c_7 (\int_{-1}^1
\parallel f(s) \parallel^p (1-s)^{\gamma p} (1+s)^{\delta p} ds
)^{1/p},
\]
and thus

\begin{eqnarray}
M_1 & \le & c_6 c_7 (\int_{-1}^1 \parallel \sum_{j=1}^{n+1} y_j \chi_j(s)
/ |I_j| \parallel^p w_{\gamma p, \delta p}(s) ds)^{1/p} \nonumber \\
 & \le & c_8 (\sum^{n+1}_{j=1} \parallel y_j \parallel^{p}/|I_j|^{p-1} w_{\gamma
 p, \delta p} (t_j))^{1/p} \\
 & \le & c_9 (\sum^{n+1}_{j=1} \lambda_j \parallel q(t_j) \parallel^p
 )^{1/p}, \nonumber
\end{eqnarray}

\noindent using that by (4) and (5)
\[
\lambda_j \sim |{p_{n+1}^{(\alpha ,\beta)}}'(t_j) |^{-p} |I_j|^{1-p}
(1-t_j)^{\gamma p} (1+t_j)^{\delta p}.
\]

The error term $M_2$ can be discretized in view of the monotonicity
properties of the $f_j$'s. The integration with respect to $t$ for $|t-t_j| \le
2|I_j|$ leads to another term $M_{21}$ of the form (13), and $M_2 \le
M_{21} + M_{22}$ with 

\begin{eqnarray}
M_{22} & = & c_{10} (\sum_{i=1}^{n+1} (\sum_{j=1, j
\ne i}^{n+1} \parallel y_j \parallel |I_j| / (t_i - t_j)^2)^p |I_i|
w_{\gamma p, \delta p} (t_i))^{1/p}\\
 & = & c_{10} [\sum_{i=1}^{n+1} (\sum_{j=1, j \ne i}^{n+1}
 a_{ij} \lambda_j^{1/p} \parallel q(t_j) \parallel )^p ]^{1/p}, \nonumber
\end{eqnarray}
 where $a_{ij} = (|I_i| / \lambda_j)^{1/p} w_{\gamma ,\delta} (t_i) |I_j|
/ (|p_{n+1}'(t_j)|(t_i - t_j)^2)$ for $i \ne j$ and $\gamma , \delta$
as before. We claim that $A_n = (a_{ij})^{n+1}_{i , j =1}$ defines a map
$A_n : \ell_p^{n+1} \longrightarrow \ell_p^{n+1}$ with norm bounded by
 a $C$ independent of $n \in \Bbb N$. Then (14) is bounded by $c_{10} C
(\sum_{j=1}^{n+1} \lambda_j \parallel q(t_j) \parallel^p)^{1/p}$ as
required. 
Calculation using (4) and (5) shows that for $i, j \le n/2$
\[
a_{ij} \sim (\frac{i}{j})^\eta \frac{j^2}{(i^2 - j^2)^2}, \qquad \qquad
\eta := (\alpha + 1/2)(2/p -1).
\]

The restriction on $\alpha $ gives that $-1/2 \le \eta \le 2$. This
easily implies \mbox{$|a_{ij}| \le c_{11} /(i-j)^{3/2}$}; for $\eta \ge 0$ or
($\eta < 0$ and $i > j/2$) one even has the bound $c_{11} / (i-j)^2$. In
any case
\[
\sup_{i \le n/2} \sum_{j \le n/2, i \ne j} |a_{ij}| \le C, \qquad
\qquad \sup_{j \le n/2} \sum_{i \le n/2, i \ne j} |a_{ij}| \le C,
\]
and hence $(a_{ij})_{i,j \le n/2}$ is uniformly bounded on
$\ell_1^{[n/2]}$ and $\ell_{\infty}^{[n/2]}$ and thus by interpolation on
$\ell_p^{[n/2]}$. The three other cases of pairs $(i,j)$, e.g. $i > n/2
\ge j$, are treated similarly, using the assumption on $\beta$ as
well.

This proves $(2) \Rightarrow (1)$ except for the case of $\alpha \le
-1/2$ or $\beta \le -1/2$ when (11) and (12) do not hold for $t \notin
J_n$. Assume e.g. $\alpha \le -1/2$. In this case, we estimate the
remaining term

\begin{eqnarray*}
M_3 &:= & (\int^1_{b_n} \parallel \sum_{j=1}^{n+1} y_j
\frac{p_{n+1}^{(\alpha ,\beta )} (t)}{t-t_j} \parallel^p w_{\alpha
\beta}(t) dt )^{1/p}
\end{eqnarray*}
by the triangle and the H\"older inequality, using (4), (5), (10) and
the fact that for $t \ge b_n$ $|t-t_j|^{-1} \le d(n/j)^2$. We find

\begin{eqnarray*}
M_3 & \le & c_{12} \ n^{-2(1+\alpha )/p} (\sum_{j=1}^{n+1} j^{\alpha -
1/2} \parallel x_j \parallel )\\
 & \le & c_{13} (\sum_{j=1}^{n+1} j^{-(\alpha (2/p -1) + 1/p +1/2)
p'})^{1/p'} (\sum_{j=1}^{n+1} \lambda_j \parallel x_j \parallel^p)^{1/p}\\
 & \le & c_{14} (\sum_{j=1}^{n+1} \lambda_j \parallel x_j \parallel^p)^{1/p}.
\end{eqnarray*}
where we have used that $\alpha (\frac{2}{p} -1)+ \frac{1}{p} - \frac{1}{2} >
 \frac{1}{p'}$.

\noindent (b) $(1) \Rightarrow (2)$. For the converse, assume the interpolation
inequality (7) to be true.

We claim that (7) implies that the Hilbert matrix $A =
((i-j+1/2)^{-1})_{i,j \in \Bbb N}$ defines a bounded operator $A :
\ell_p (X) \to \ell_p (X)$. A well-known approximation and scaling
argument shows that this is equivalent to the boundedness of the Hilbert
transform $H$ in $L_p(\Bbb R; X)$, i.e. $X$ is a UMD-space and
necessarily $1 < p < \infty$. In this sense $A$ is a discrete version
of $H$. For $n \in \Bbb N$ we need the zeros $(t_j^{n+1})_{j=1}^{n+1}$
of $p_{n+1}^{(\alpha ,\beta)}$ and $(t_i^n)^n_{i=1}$ of $p_n^{(\alpha
,\beta)}$, ordered decreasingly as before and the corresponding quadrature
weights $\lambda_j^{n+1}$ and $\lambda^n_i$. Let \mbox{$J_n:=\{ j \in \Bbb N \ | 
\ n/4
\le j \le 3n/4\}$}. For any sequence $(x_j)_{j \in J_n} \subseteq X$,
consider the $X$-valued polynomial
\[
q:= \sum_{j \in J_n} \lambda_j^{-1/p} x_j \ell_j^{n+1} \in \Pi_n
(X) \qquad ,\qquad \ell_j^{n+1} (t_i^{n+1}) = \delta_{ij}.
\]
Applying (7) to $q$ and inequality (6) with $(n+1)$ replaced by $n$ we
find
\[
(\sum_{i \in J_n} \lambda^n_i \parallel q (t_i^n) \parallel^p )^{1/p}
\le c_1 (\int_{-1}^1 \parallel q(t) \parallel^p w_{\alpha \beta} (t) dt
)^{1/p} \le c_2 (\sum_{j \in J_n} \lambda_j^{n+1} \parallel q (t_j^n)
\parallel^p )^{1/p},
\]
i.e.
\begin{equation}
(\sum_{i \in J_n} \parallel \sum_{j \in J_n} b_{ij} x_j \parallel^p
)^{1/p} \le c_2 (\sum_{j \in J_n} \parallel x_j \parallel^p )^{1/p}
\end{equation}
with
\[
b_{ij} := (\lambda_i^n / \lambda_j^{n+1} )^{1/p} \ell_j^{n+1} (t_i^n) =
(\lambda_i^n / \lambda_j^{n+1} )^{1/p} \frac{p_{n+1}^{(\alpha ,\beta)}
(t_i^n)}{{p_{n+1}^{(\alpha ,\beta )}}'(t_j^{n+1})(t_j^n -
t_j^{n+1})}.
\]

Let $k_n := |J_n | \sim n/2$. By (15), $B_n := (b_{ij})_{i,j \in J_n}$
satisfies $\parallel B_n : \ell_p^{k_n}(X) \longrightarrow \ell_p^{k_n}
(X) \parallel \le c_2$, $c_2$ being independent of $n \in \Bbb N$. $B_n$
is close to the block $A_n := ((i-j+1/2)^{-1})_{i,j \in J_n}$ of the
Hilbert matrix $A$. To show this we first evaluate $b_{ij}$. By
Szeg\"o \cite[(4,5.7)]{kn:sz}
\[
(1-t^2){p_n^{(\alpha ,\beta )}}' (t) = ({\eta_n}' t + {\eta_n}'')(
p_n^{(\alpha ,\beta )} (t) - \eta_n p_{n+1}^{(\alpha ,\beta )} (t),
\]
where $\eta_n , {\eta_n}' , {\eta_n}'' \in \Bbb R$ depend on $n$ and
$(\alpha ,\beta )$, with $\eta_n /n \longrightarrow 1$ for $n
\longrightarrow \infty$.
Hence, using (5), for $i \le 3n/4$

\begin{equation}
p_{n+1}^{(\alpha ,\beta )} (t_i^n) = \eta^{-1}_n
(1-(t_i^n)^2){p_n^{(\alpha ,\beta )}}' (t_i^n) \ \sim \ (-1)^i
(n/i)^{\alpha - 1/2}.
\end{equation}

Thus $b_{ij} = \gamma_i^n/(\delta_j^n n(t_i^n -
t_j^{n+1}))$ where $0 < c_3 \le |\gamma^n_i |, |\delta^n_j | \le c_4 <
\infty$ for $i,j \in J_n$ and the matrices $C_n:= (c_{ij})_{i,j \in
J_n}$, $c_{ij} := n^{-1}(t_i^n - t_j^{n+1})^{-1}$, are uniformly bounded
on $\ell_p^{k_n}(X)$ as well. By Szeg\"{o} \cite[(8.9.8]{kn:sz}, 

\begin{equation}
t_i^n = \cos \theta^n_i, \qquad \qquad \theta_i^n =
\frac{(i+k+\varepsilon_{ni}) \pi}{n+\alpha + 1/2},
\end{equation}

\noindent where $k$ depends on $(\alpha ,\beta)$ only and $\sup_{i \in J_n} 
|\varepsilon_{ni}|
\longrightarrow 0$ for $n \longrightarrow \infty$. For $k_n \times k_n$
matrices $D_n$ and $E_n$, we write $D_n \approx E_n$ provided that the
matices $D_n - E_n$ are uniformly bounded as maps on $\ell_p^{k_n}(X)$,
for {\em any} $1 \le p \le \infty$. It suffices to show that $C_n
\approx A_n$, then $\sup_{n \in \Bbb N} \parallel A_n : \ell_p^{k_n} (X)
\longrightarrow \ell_p^{k_n} (X) \parallel \le c_5$, i.e. $X$ is a
UMD-space. We have
\[
c_{ij} = \frac{1}{2n} \frac{1}{\sin (\theta_i^n - \theta_j^{n+1})/2}
 \ \frac{1}{\sin (\theta_i^n + \theta_j^{n+1})/2}.
\]
If
\[
d_{ij} := \frac{1}{\pi (i-j+1/2) \sin (\theta_i^n
+\theta_j^{n+1})/2}, \qquad D_n := (d_{ij})_{i,j \in J_n},
\]

\noindent then $C_n \approx D_n : 1/\sin [(\theta_i^n + \theta_j^{n+1})/2]$ is
uniformly bounded for $i,j \in J_n$, and hence the estimates $|1/\sin x -
1/x | \le |x|$ for $|x| \le \pi /4$ and  $|1/x - 1/(x+\varepsilon
)| \le 2\varepsilon/x^2$ for $x = i-j + 1/2$ and $|\varepsilon | \le
1/4$ yield

\begin{eqnarray*}
|c_{ij} - d_{ij} | & \le & c_6 (| \frac{1}{2n \sin (\theta_i^n -
\theta_j^{n+1})/2} - \frac{1}{n (\theta_j^n - \theta_j^{n+1})} | + |
\frac{1}{n(\theta_i^n - \theta_j^{n+1})} - \frac{1}{\pi(i-j+1/2)} | ) \\
 & \le & c_7 (\frac{|i-j+1/2|}{n^2} + \frac{1}{|i-j+1/2|^2}).
\end{eqnarray*}

Hence $\sup_i \sum_j |c_{ij} - d_{ij}| \le c_8 , \ \sup_j \sum_i |c_{ij}
- d_{ij} | \le c_8$ uniformly in $n \in \Bbb N$, i.e. $C_n \approx D_n$
(first for $p=1, \infty$, then by interpolation for general $p$). Next
$D_n$ is transformed into $E_n = (e_{ij})$ with $D_n \approx E_n$,

\[
e_{ij} := \frac{1}{\pi (i-j+1/2) \sin \frac{i+j}{2n} \pi};
\qquad  i,j \in J_n.
\]

By the Lipschitz continuity of $1/\sin x$ in $\pi /4 \le x \le 3 \pi
/4$,

\[
|d_{ij} - e_{ij} | \le c_9 / (n|i-j+1/2|),
\]

\noindent which again is uniformly row- and column- summable, i.e. $D_n \approx
E_n$. Finally, let

\[
f_{ij} := \frac{1}{\pi (i-j+1/2) \sin \frac{i}{n} \pi}, \qquad
 \qquad F_n := (f_{ij})_{i,j \in J_n}.
\]

For $y \in [\pi /8, 3\pi /8]$, $g(x) = 1/\sin (x+y)$ is
Lipschitz-continous in $x \in [\pi /8, 3\pi /8]$ with constant $\le 2$.
Hence for $i,j \in J_n$

\[
|e_{ij} - f_{ij} | \le \frac{|i-j|}{n |i-j+1/2|} \le
\frac{2}{n} , \qquad E_n \approx F_n.
\]

Since $(\sin (\frac{i \pi}{n})^{-1})_{i \in J_n}$ is bounded away from
zero, this implies that \mbox{$A_n = ((i-j+1/2)^{-1})_{i,j \in J_n}$} defines
uniformly bounded maps on $\ell_p^{k_n} (X)$. Hence $X$ is a UMD-space
and $1 < p < \infty$.

We now prove that (7) implies $\mu (\alpha ,\beta ) < p < M(\alpha
,\beta )$. This is a purely scalar argument, $X = \Bbb R$. By symmetry,
we may assume that $\alpha \ge \beta$; the case of $\alpha \le -1/2 \ \
(1<p<\infty )$ is known already. So let $-1/2 < \alpha$.
Take $q = p_n^{(\alpha ,\beta )} \in \Pi_n$ in (7) to show that
necessarily $p < M(\alpha ,\beta ) = 4 \frac{\alpha +1}{2\alpha +1}$:
Using the second formula of Szeg\"{o} \cite[4.5.7]{kn:sz}, one shows
similarly as in (16) that $|p_n^{(\alpha ,\beta )} (t_i^{n+1}) | \sim
(n/i)^{\alpha -1/2}$ for $i \le n/2$ and $\sim (n/(n+2-i))^{\beta -1/2}$
for $i > n/2$.
Thus by (4) and Newman-Rudin \cite{kn:nr}, cf. also (10),

\begin{equation}
(\int_{-1}^1 |p_n^{(\alpha ,\beta )} (t) |^p w_{\alpha \beta} (t)
dt)^{1/p} \sim \left\{\begin{array}{cl}
1 & p < M (\alpha ,\beta ) \\
(\log n)^{1/p} & p=M(\alpha ,\beta ) \\
n^{p(\alpha + 1/2)-2(\alpha +1)} & p>M(\alpha ,\beta )
\end{array}\right\},
\end{equation}

\begin{eqnarray}
\lefteqn{(\sum_{j+1}^{n+1} \lambda_j |p_n^{(\alpha ,\beta )} (t_j) |^p )^{1/p} 
\sim }\\
 & \left\{\begin{array}{cl}
1 & \mbox{($p<\infty$ and $\alpha \le 1/2$) or ($p < 4 (\alpha
+1)/(2\alpha -1)$ and $\alpha > 1/2$)}\\
(\log n)^{1/p} & \mbox{$p=4(\alpha+1)/(2\alpha -1)$ and $\alpha > 1/2$}\\
n^{(\alpha - 1/2)-2/p(\alpha +1)} & \mbox{$p > 4(\alpha +1)/(2\alpha -1)$
 and $\alpha > 1/2$}
\end{array}\right\}. & \mbox{} \nonumber
\end{eqnarray}

Hence for $p \ge M(\alpha ,\beta )$, the order of growth (in $n$) in
(18) is faster than in (19) and (7) cannot hold. To prove that
necessarily $p > \mu (\alpha ,\beta ) = 4(\alpha +1) / (2\alpha + 5)$
for $\alpha > 1/2$ ($\mu (\alpha , \beta ) = 1$ for $\alpha \le 1/2$), we
take
\[
q= \ell_1 \in \Pi_n , \qquad \ell_1 (t) = p_{n + 1}^{(\alpha ,\beta )}
(t) / ((t-t_1){p_{n+1}^{(\alpha ,\beta )}}'(t_1)).
\]
Clearly, the right side of (7) is $\sim \lambda_1^{1/p} \sim
n^{-2/p(\alpha +1)}$ by (4) whereas the asymptotic formulas for
$p_n^{(\alpha ,\beta)}$ of Szeg\"o \cite[8.21]{kn:sz} and (5) yield

\begin{eqnarray}
(\int^1_{-1} |\ell_1 (t) |^p w_{\alpha \beta} (t) dt )^{1/p} & \sim & 
(\int_0^{1-n^{-2}} (1-t)^{\alpha - p/4(2 \alpha +5)} dt)^{1/p} 
n^{-(\alpha + 5/2)} \nonumber \\
 & \sim & \left\{\begin{array}{cl}
n^{-(\alpha + 5/2)} & p < \mu (\alpha ,\beta ) \\
(\log n)^{1/p} & p= \mu (\alpha ,\beta ) \\
n^{-2/p (\alpha ,\beta )} & p > \mu(\alpha,\beta)
\end{array}\right\}.
\end{eqnarray}

Hence (20) grows faster in $n$ than $\lambda_1^{1/p} \sim n^{-2/p(\alpha
+ 1)}$ if $p \le \mu (\alpha ,\beta )$, i.e. $p > \mu (\alpha ,\beta )$
is necessary for (7) to hold. This proves (b) of Theorem 1.

(a). The left interpolation inequality (6) in Theorem 1 is proved
as in the scalar case. Nevai's proof in \cite{kn:n} using the mean value
theorem, H\"older's inequality and some weighted form of Bernstein's
inequality in the p-norm generalizes directly to the vector-valued
setting. Just as the scalar result of Khalilova \cite{kn:kl} and Potapov
\cite{kn:pot}, the vector-valued form of the Bernstein $L_p$-inequality
(lemma 2 in \cite{kn:n}) is proved by interpolating at Tchebychev nodes,
using an averaging technique, the triangle inequality in $L_p$ and the
Bernstein inequality for the sup-norm. In the vector valued case the
latter follows from the scalar version, applying linear functionals and
using the Hahn-Banach theorem. We do not give the details, since the proofs 
of \cite{kn:n}, \cite{kn:kl} and \cite{kn:pot} directly generalize.
{\eproof}
\noindent {\bf Remarks.}

(1). If the validity of (7) of Theorem 1 (b),
\[
(\int_I \parallel q (t) \parallel^p w_{\alpha \beta} (t) dt)^{1/2} \le
c_p (\sum_{j=1}^{n+1} \lambda_j^{n+1} \parallel q (t_j^{n+1}) \parallel^p
)^{1/p},
\]
is required only for all polynomials $q \in \Pi_k (X)$ with $k \le n/2$,
this holds for all $1 \le p \le \infty$ and all Banach spaces, at least
if $\alpha ,\beta \ge -1/2$. This follows from the boundedness of the
generalized de la Vall\'ee-Poussain means in $L_p (X)$ along similar
lines as in Zygmund \cite{kn:z2}, Stein \cite{kn:st} and Askey
\cite{kn:as}. Thus the restriction on $p$ and $X$ in Theorem 1 comes
from requiring the number of nodes to equal the dimension of $\Pi_n$,
namely $(n+1)$. In this way, however, one isomorphically identifies
$\Pi_n (X) \subset L_p (X)$ with the space $\ell_p^{n+1} (X)$.

(2). The proof of the necessity of the UMD-condition for
inequality (7) of Theorem 1 (b) will work for more general orthogonal
polynomials provided that sufficiently precise information on a fairly
large part of the zeros of these is known, like in (17).

The restriction $p < M(\alpha ,\beta )$ means geometrically (for $\alpha
> -1/2$) that the value $|p_n^{(\alpha ,\beta )} (t_1^{n+1})|^p $ is
much smaller than the dominating mean value of $|p_n^{(\alpha ,\beta
)} (t) |^p$ over $(t_1^{n+1} ,1)$ with respect to $w_{\alpha
\beta} (t) dt$, if $p \ge M(\alpha ,\beta)$.

An immediate corollary to Theorem 1 is the following result on the
convergence of interpolating polynomials which in the scalar case is due
to Askey \cite{kn:as} and Nevai \cite{kn:n}.

\begin{proposition}
Let $X$ be a UMD-space, $\alpha ,\beta > -1$ and $p < M(\alpha ,\beta
)$. Let \mbox{$f: (-1,1) \longrightarrow X$} be continuous. Then the
interpolating polynomials of $f$ at the zeros $(t_j)_1^{n+1}$ of
$p_{n+1}^{(\alpha ,\beta )}$, \mbox{$I_n f: = \sum_{j=1}^{n+1} f(t_j) \ell_j \in
\Pi_n (X)$}, converge to $f$ in the p-norm,
\[
\parallel f - I_n f \parallel_{p;\alpha ,\beta } = (\int_{-1}^1
\parallel f(t) - I_n f(t) \parallel^p w_{\alpha \beta} (t) dt)^{1/p}
\longrightarrow 0.
\]
\end{proposition} 
{\kproof} Approximate $f$ by polynomials $q_n \in \Pi_n
(X)$ in the sup-norm, $\parallel f-q_n \parallel_{\infty}
\longrightarrow 0$. We may assume that $\mu (\alpha ,\beta ) < p <
M(\alpha ,\beta )$, since the p-norms get weaker for smaller $p$. Using
(b) of Theorem 1 and
\[
\sum_{j=1}^{n+1} \lambda_j = \int_{-1}^1 w_{\alpha \beta} (t) dt =: M <
\infty,
\]
we find
{\samepage \begin{eqnarray*}
\parallel f -I_n f \parallel_{p; \alpha ,\beta } & \le & \parallel f-q_n
\parallel_{p ;\alpha ,\beta} + \parallel q_n -I_n f \parallel_{p:\alpha
,\beta }\\ 
 & \le & M^{1/p} \parallel f-q_n \parallel_{\infty} + c_p
(\sum_{j=1}^{n+1} \lambda_j \parallel q_n (t_j) - f(t_j) \parallel^p
)^{1/p} \\ & \le & (1 +c_p) M^{1/p} \parallel f-q_n \parallel_{\infty}
\longrightarrow 0.
\end{eqnarray*}
{\eproof}}

\section{Convergence of vector-valued Jacobi series}
{\bproof} {\bf of Theorem 2:} Recall that $Q_nf := \sum_{j=0}^n <
f,p_j^{(\alpha ,\beta )} > p_j^{(\alpha ,\beta )}$ for $f \in L_p (I,
w_{\alpha ,\beta} ; X)$. Thus $Q_n$ is the integral operator induced by
the kernel $k_n(x,y)=\sum^n_{j=0} p_j^{(\alpha,\beta)}(x) p_j^{(\alpha ,\beta 
)} (y)$ with
respect to the measure $d\mu (t) = w_{\alpha \beta } (t) dt.$

\noindent (2) $\Rightarrow$ (1). We sketch the straightforward generalization of 
the
scalar proof of Pollard \cite{kn:pol} and Muckenhaupt \cite{kn:mu} to
the UMD-case. Since $Q_nf \rightarrow f$ on the dense set of $X$-valued
polynomials $f$, (1) of Theorem 2 is equivalent to
\begin{equation}
\sup_{n \in \Bbb N} \parallel Q_n:L_p (I, w_{\alpha ,\beta }; X)
\longrightarrow L_p (I, w_{\alpha \beta} , X ) \parallel = c_p < \infty.
\end{equation}

Using the Christoffel-Darboux formula for $k_n$ and the classical
analysis of Pollard \cite{kn:pol}, (21) will follow from the uniform
boundedness of the integral operators $T_{n1}$, $T_{n2}$, $T_{n3}$ induced by
the following kernels as maps in $L_p (I,w_{\alpha \beta }; X)$:
\[
k_{n1} (x,y) := p_{n+1}^{(\alpha ,\beta )} (x) q_n^{(\alpha ,\beta )}
(y)/(x-y) \ \ , \ \ q_n^{(\alpha ,\beta )} (y) := (1-y^2) p_n^{(\alpha
+1,\beta +1)} (y)
\]
\[
k_{n2} (x,y) := k_{n1}(y,x) \ \ , \ \ k_{n3}(x,y) := p_n^{(\alpha ,\beta
)} (x) p_n^{(\alpha ,\beta )} (y).
\]

The proof of the uniform boundedness of $T_{n1}$ and $T_{n2}$ is similar to
the proof of Theorem 1, (b), (2) $\Rightarrow$ (1). On the intervals
$J_n$ defined there (for $\alpha ,\beta \ge -1/2$, $J_n =I$), $T_{n1}$ and
$T_{n2}$ are uniformly bounded in p-norm provided that the weighted
Hilbert transform kernels
\mbox{$(x-y)^{-1} (w_{\alpha \beta} (x) w_{\alpha \beta} (y))^{1/2}
((1-x^2)/(1-y^2))^{\pm 1/4}$} (+ for $T_{n2}$, - for $T_{n1}$) define
bounded operators on $L_p (I;X)$, as follows from (10) the same way as in
(b), (2) $\Rightarrow$ (1). In view of the UMD-assumption on $X$, this
will follow from the boundedness of the kernel operator defined by
\[
\frac{1}{|x-y|} |(w_{\alpha \beta} (x) w_{\alpha \beta } (y))^{1/2}
((1-x^2)/(1-y^2))^{\pm1/4} -1|
\]
on $L_p(I;X)$. Using again lemma 1, the latter fact is a consequence of 
\[
-\frac{1}{p} < \alpha (\frac{1}{p} - \frac{1}{2} ) \pm \frac{1}{4} ,
\beta (\frac{1}{p} - \frac{1}{2}) \pm \frac{1}{4} < 1 - \frac{1}{p},
\]
i.e. $m(\alpha ,\beta) < p < M(\alpha ,\beta )$. If e.g. $\alpha <
-1/2$, the part of $\parallel T_{ni}f \parallel_p$, $i \in \{ 1,2 \}$,
on the interval $(1-n^{-2},1)$ outside $J_n$ has to be estimated
separately. However, $p_n^{(\alpha ,\beta)}$ and $q_n^{(\alpha ,\beta )}$
are uniformly bounded in $n \in \Bbb N$ there, and a direct application
of the continuity of the (unweighted) Hilbert transform suffices. The
uniform boundedness of $T_{n3}$ follows from $\sup_{n \in \Bbb N}
\parallel p_n^{(\alpha ,\beta )} \parallel_p \parallel p_n^{(\alpha
,\beta )}  \parallel_{p'} < \infty$ if $m(\alpha ,\beta ) < p <
M(\alpha ,\beta )$.

(1) $\Rightarrow $ (2). Assume that $Q_n f \longrightarrow f$ for all $f
\in L_p(I, w_{\alpha \beta} ;X)$. By the Banach-Steinhaus theorem, this
is equivalent to (21). Using (21), we prove (7) of Theorem 1, which 
then implies that $X$ is an UMD-space and, in view of the self-duality
of (21), that $m(\alpha ,\beta ) < p < M(\alpha ,\beta )$. To show (7),
we dualize (6) which holds for all $X$ and $p$. Let $q \in \Pi_n (X)$.
Then there is a $g \in L_{p'} (I, w_{\alpha \beta }; X^{\ast} )$ which
$\parallel g \parallel_{p';\alpha ,\beta} =1$ and 

\begin{eqnarray*}
J:=(\int^1_{-1} \parallel q(t) \parallel^p w_{\alpha \beta} (t) dt
)^{1/p} & = & \int^1_{-1} <q(t),g(t) >_{(X,X^{\ast})} w_{\alpha \beta} (t)
dt \\
  & = & \int^1_{-1} <q(t), Q_ng(t) >_{(X,X^{\ast})} w_{\alpha \beta} (t)
 dt.
\end{eqnarray*}

Since $<q,Q_ng> \in \Pi_{2n}$, Gaussian quadrature, H\"older's
inequality and (6) as well as (the dual form of) (21) yield

\begin{eqnarray*}
J & = & \sum_{j=1}^{n+1} \lambda_j <q(t_j), Q_ng(t_j)> \\
 & \le & (\sum_{j=1}^{n+1} \lambda_j \parallel Q_n g(t_j)
 \parallel^{p'}_{X^{\ast}} )^{1/p'} (\sum_{j=1}^{n+1} \lambda_j
 \parallel q(t_j) \parallel_X^p )^{1/p} \\
 & \le & c \parallel Q_ng \parallel_{p';\alpha ,\beta } (\sum_{j=1}^{n+1}
 \lambda_j \parallel q(t_j) \parallel^p )^{1/p} \\
 & \le & c \, c_p (\sum_{j=1}^{n+1} \lambda_j \parallel q(t_j) \parallel^p
 )^{1/p} 
\end{eqnarray*}
which is (7).
{\eproof}

We turn to the equivalence of vector-valued Jacobi means.

{\bproof} {\bf of proposition 3:} For $\alpha > -1$ and $w_{\alpha} =
w_{\alpha ,\alpha} $, the map \mbox{$\psi :L_2 (I, w_{\alpha }; X)
\longrightarrow L_2 (0, \pi ; X)$} defined by 
\[
\psi (g)(s) = (\sin s )^{\alpha + 1/2 } g (\cos s) \qquad g \in L_2(I,w_{\alpha
};X), \qquad s \in [0,\pi],
\]
is an isometry. Let $q_n^{(\alpha )} := \psi (p_n^{(\alpha ,\alpha )} )$
and
\[
{\cal K}_n^{(\alpha ,\beta )} (t,s) := \sum_{j=0}^n q_j^{(\alpha )} (t)
q_j^{(\beta )} (s) ; \qquad \alpha, \beta > -1.
\]

These kernels induce uniformly bounded operators on $L_p (0, \pi ;X)$
for any $1<p<\infty$, e.g. there is $c_p$ such that for all $n \in \Bbb
N$ and $h \in L_p (0, \pi ;X)$
\begin{equation}
(\int_0^{\pi} \parallel \int_0^{\pi} {\cal K}_n^{(\alpha ,\beta )} (t,s) h(s)
ds \parallel^p dt )^{1/p} \le c_p (\int_0^{\pi} \parallel h(s)
\parallel^p ds )^{1/p}.
\end{equation}

This follows from the proofs of Theorem 1 and 3 of Gilbert \cite{kn:gi}:
The scalar proof given there directly generalizes to the $X$-valued
UMD-case since only the $L_p$-uniform boundedness of the Dirichlet and
conjugate Dirichlet kernel operators is used, which holds $X$-valued for
UMD-spaces. In effect, ${\cal K}_n^{(\alpha ,\beta )}$ is shown in \cite{kn:gi}
to behave very similar to the Dirichlet kernel. In particular
\begin{equation}
|{\cal K}_n^{(\alpha ,\beta )} (t,s) | \le d_1 /|t-s|
\end{equation}
where $d_1$ is independent of $n \in \Bbb N$ and $t, s \in [0, \pi ]$.
We claim that also 
\begin{equation}
(\int_0^{\pi} \parallel \int_0^{\pi} {\cal K}_n^{(\alpha ,\beta )} (t,s) (\sin
t/ \sin s)^{\gamma + 1/p -1/2} h(s) ds \parallel^p dt )^{1/p} \le {c_p}'
(\int_0^{\pi} \parallel h(s) \parallel^p ds )^{1/p},
\end{equation}
provided that $|\gamma /2 + 1/p - 1/2 | < 1/4$. By (22), this will
follow from the uniform boundedness of the difference kernel operators
\[
{\cal L}_n^{(\alpha ,\beta )} (t,s) := {\cal K}_n^{(\alpha ,\beta )} (t,s) ((\sin t
/\sin s)^{\gamma + 1/p -1/2} -1)
\]
in $L_p (0,\pi ;X)$. Using (23) and elementary estimates we obtain the 
existence of a $d_2$ such that for $n \in \Bbb N$ and $t,s
\in [0, \pi /2]$
\[
|{\cal L}_n^{(\alpha ,\beta )} (t,s) | \le \frac{d_1}{|t-s|} | (\frac{\sin
t}{\sin s})^{\gamma + 1/p - 1/2} -1 | \le \frac{d_2}{|t-s|}
|(\frac{t}{s})^{\gamma + 1/p - 1/2 } -1 |.
\]
Hence by lemma 1, the ${\cal L}_n^{(\alpha
,\beta )}$-kernels define uniformly bounded integral operators in
$L_p(0, \pi/2 ; X)$ since $-1/p < \gamma + 1/p - 1/2 < 1 - 1/p$. On
$(\pi /2, \pi )$, the estimate is similar; for $t \in [\pi /2 ,\pi ], s
\in [0, \pi /2 ]$, there are only point singularities and the
transformation $t \rightarrow \pi - t$ reduces the ${\cal L}_n^{(\alpha
,\beta )}$-boundedness to the one of the positive kernel $1/(t+s)$ in
$L_p (0, \pi /2 ;X)$. Hence (24) holds.

For functions $f \in L_p (I,w_{(\beta + \gamma ) p/2} ; X)$ on the interval
$I = (-1,1)$ and the kernel
\[
k_n^{(\alpha ,\beta )} (x,y):= \sum^n_{j=0} p_j^{(\alpha ,\alpha )} (x)
p_j^{(\beta ,\beta )} (y),
\]
(24) is equivalent to 
\begin{eqnarray}
( \int^1_{-1} \parallel \int^1_{-1} k_n^{(\alpha ,\beta )} (x,y) f(y)
w_{\beta} (y) dy \parallel^p w_{(\alpha + \gamma )p/2} (x) dx )^{1/p}
\nonumber \\ 
\le {c_p}' \int^1_{-1} \parallel f(y) \parallel^p w_{(\beta + \gamma
)p/2} (y) dy )^{1/p} 
\end{eqnarray}
as the transformation $h(s) = (\sin s)^{\beta + \gamma + 1/p} f(\cos s)$
shows. Applying (25) to $f(y)= \sum^n_{j=0} p_j^{(\beta ,\beta )}(y)
x_j$, where $x_j \in X$, yields a one-sided estimate of (8); the
converse direction follows from the symmetry of the statement in $\alpha
$ and $\beta $. The argument also shows that the convergence of the
series $\sum_j p_j^{(\alpha ,\alpha )} \otimes x_j$ in $L_p (I,
w_{(\alpha + \gamma )p/2} ; X)$ is equivalent to the convergence of the series
$\sum_j p_j^{(\beta ,\beta )} \otimes x_j $ in $L_p (I, w_{(\beta +
\gamma )p/2} ; X)$, provided that \mbox{$|\gamma/2 + 1/p - 1/2 | < 1/4$}.
{\eproof}

The choice of $p=2$ and $\gamma =0$ shows that the means 
\[
(\int^1_{-1} \parallel \sum_{j=0}^n p_j^{(\alpha ,\alpha )} (t) x_j
\parallel^2 w_{\alpha } (t) dt )^{1/2}
\]
are essentially independent of $\alpha $, the choice of $\gamma =0$ for
$4/3 < p < 4$ shows a similar statement for the means
\[
(\int^1_{-1} \parallel \sum^n_{j=0} p_j^{(\alpha ,\alpha )} (t) x_j
\parallel^p w_{\alpha p/2} (t) dt)^{1/p}.
\]

\section{Unconditional convergence}
We now show that under the conditions of Proposition 4, vector-valued
convergence of orthonormal series is unconditional only in the case of
Hilbert spaces.

{\bproof} {\bf of Proposition 4:}

(i). Let $(\Omega , \mu )$ be a finite measure space
and $(p_n)$ be a complete orthonormal system in $L_2(\Omega ,\mu )$ such
that $\sum_n < f,p_n > p_n$ converges unconditionally for all $f \in L_p
(\Omega ,\mu )$. By duality, the same holds in $L_{p'}(\Omega ,\mu )$.
Thus we may assume that $p \ge 2$. Using the unconditionality and the
Khintchine inequality, we find for any finite sequence $(a_n) \subset
\Bbb K$
\begin{eqnarray*}
\sum_n |a_n |^2 = \parallel \sum_n a_n p_n \parallel^2_2 & \le & c_1
\parallel \sum_n a_n p_n \parallel^2_p \\ 
  & \le & c_2 (\int^1_0 \int_{\Omega} |
\sum_n a_n r_n (t) p_n (w) |^p d\mu (w) dt )^{2/p} \\
 & \le & c_3 (\int_{\Omega}
(\sum_n |a_n |^2 |p_n (w) |^2 )^{p/2} d\mu (w))^{2/p} \\
 & \le & c_3 \sum_n
(\int_{\Omega} |a_n |^p |p_n (w) |^p d\mu (w))^{2/p} \\
 & = & c_3 \sum_n |a_n
|^2 \parallel p_n \parallel^2_p \ \le \ c_3 \sup_n \parallel p_n
\parallel^2_p \sum_n |a_n |^2 \\
 & \le & c_4 \sum_n |a_n|^2,
\end{eqnarray*}
where we have used the triangle inequality in $L_{p/2}$ and the
assumption on $\parallel p_n \parallel_p$.

Hence $\parallel\sum_n a_n p_n \parallel_p \sim \parallel (a_n) \parallel_{l_2}$,
which implies $p=2$ since $(p_n)$ was assumed to be a {\em complete}
orthonormal system and $L_p(\Omega ,\mu ) \sim \ell_2$ only for $p=2$.
Thus (i) holds, even for $X = \Bbb K$.

(ii). We give a modification of the argument of Defant and Junge
\cite{kn:dj}. Let \mbox{$x_1,\cdots x_m \in X$}. By the unconditionality
assumption on the $(p_j)$, the hypothesis that \mbox{$\sup_j | p_j | \in L_2
(\Omega ,\mu )$}, and the contraction principle, cf. Maurey and Pisier
\cite{kn:mp} we get,

\begin{eqnarray*}
(\int_{\Omega} \parallel \sum^m_{j=1} p_j(w)x_j \parallel^2 d\mu (w)
)^{1/2} & \le & c_1 (\int_{\Omega } \int^1_0 \parallel \sum^{m}_{j=1}
r_j (t) p_j (w) x_j \parallel^2 dt d\mu (w))^{1/2} \\
 & \le & c_2 ( \int_{\Omega } (\sup_j |p_j (w) |^2)(\int_0^1 \parallel
 \sum^{m}_{j=1} r_j (t) x_j \parallel^2 dt) d\mu (w))^{1/2} \\
 & \le & c_3 (\int^1_0 \parallel \sum^{m}_{j=1} r_j (t) x_j \parallel^2
 dt )^{1/2}.
\end{eqnarray*}

Let $(\gamma_j )$ be a sequence of independent standard
$N(0,1)$ Gaussian variables on a probability space $(\Gamma, \nu)$. By
Pisier \cite{kn:pi1} with $c_4=\sqrt{\pi/2} \, c_3$
\begin{equation}
(\int_{\Omega } \parallel \sum^{m}_{j=1} p_j (w) x_j \parallel^2 d\mu
(w))^{1/2} \le c_4 (\int_{\Gamma } \parallel \sum^m_{j=1} \gamma_j (s)
x_j \parallel^2 d\nu (s))^{1/2}
\end{equation}

Since $L_2 (\Omega ,\mu )$ is infinite dimensional, for any $n \in \Bbb
N$ there is a unitary map \\ 
\noindent $U: L_2 (\Omega ,\mu ) \longrightarrow L_2
(\Omega ,\mu )$ such that with $f_j := Up_j$ the functions $f_1, \cdots
f_n$ are mutually disjointly supported. Let $f_k = \sum_j u_{jk} p_j \in
L_2 (\Omega ,\mu )$, $(u_{jk})$ unitary. Applying (26) for arbitrary 
$y_1, \cdots y_n\in X$ with $x_j := \sum^n_{k=1} u_{jk} y_k$ we find using
the unitary invariance of the right side of (26)

\begin{eqnarray*}
(\sum^n_{k=1} \parallel y_k \parallel^2 )^{1/2} & = & (\int_{\Omega }
\parallel \sum^n_{k=1} f_k(w) y_k \parallel^2 d\mu (w))^{1/2} \\
 & = & (\int_{\Omega } \parallel \sum_j p_j(w) x_j \parallel^2 d\mu
 (w))^{1/2} \\
 & \le & c_4 (\int_{\Gamma } \parallel \sum_j \gamma_j (s) x_j \parallel^2
 d\nu (s))^{1/2} \\
 & = & c_4 (\int_{\Gamma } \parallel \sum^n_{k=1} \gamma_k (s) y_k
 \parallel^2 d\nu (s))^{1/2},
\end{eqnarray*}
i.e. $X$ has cotype 2. Similarly, the converse inequality to (26) will
imply that $X$ has type 2 and thus by Kwapie\'n \cite{kn:kw} that $X$ is
isomorphic to a Hilbert space. By Maurey and Pisier \cite{kn:mp}, the
Gaussian and the Rademacher means are equivalent since $X$ has cotype 2.
Using this and Kahane's inequality \cite{kn:lt}, we get for any $x_1 ,
\cdots x_m \in X$

\begin{eqnarray*}
(\int_{\Gamma } \parallel \sum^m_{j=1} \gamma_j (s) x_j \parallel^2
d\nu (s))^{1/2} & \le & c_5 (\int^1_0 \parallel \sum^m_{j=1} r_j (t) x_j
\parallel^2 dt )^{1/2} \\
 & \le & c_6 \int_0 ^1 \parallel \sum^m_{j=1} r_j (t) x_j \parallel dt.
\end{eqnarray*}

Since $\parallel p_j \parallel_2 =1$, the contraction principle, the
H\"older inequality and the unconditionality assumption yield similarly
as in Defant and Junge \cite{kn:dj}, cf. also Pisier \cite{kn:pi1}, that this
is

\begin{eqnarray*}
 & = & c_6 \int^1_0 \parallel \sum^m_{j=1} (\int_{\Omega } |
p_j (w) |^2 d\mu (w)) r_j (t) x_j \parallel dt \\
 & \le & c_7 \int^1_0 \int_{\Omega } (\sup_j |p_j (w) |) \parallel
 \sum^m_{j=1} p_j (w) r_j (t) x_j \parallel d\mu (w) dt \\
 & \le & c_7 \parallel \sup_j |p_j | \parallel_{L_2(\Omega ,\mu )}
 (\int_{\Omega } \int^1_0 \parallel \sum^m_{j=1} r_j (t) p_j (w) x_j
 \parallel^2 dt d\mu (w))^{1/2} \\
 & \le & c_8 (\int_{\Omega } \parallel \sum^m_{j=1} p_j (w) x_j \parallel^2
 d\mu (w))^{1/2},
\end{eqnarray*}
i.e. the converse to (26) holds. Hence $X$ is isomorphic to a Hilbert
space.
{\eproof}

As a corollary we find
\begin{proposition}
Let $\alpha ,\beta > -1$, $1 \le p \le \infty $ and $X$ be a Banach space.
Assume that for all $f \in L_p (I, w_{\alpha \beta } ; X)$, the Jacobi
series $\sum^{\infty }_{n=0} < f,p_n^{(\alpha ,\beta )} > p_n^{(\alpha
,\beta )}$ converges unconditionally in $L_p (I, w_{\alpha \beta } ;
X)$. Then $p=2$ and $X$ is isomorphic to a Hilbert space: the expansions
converge unconditionally precisely in the Hilbert space situation.
\end{proposition}
 
{\kproof} 
By Theorem 2, necessary for convergence is $m(\alpha ,\beta ) < p <
M(\alpha ,\beta )$. For these values of $p$, the inequality (10) yields
e.g. if $p>2$

\begin{eqnarray*}
\parallel p_j^{(\alpha ,\beta )} \parallel_{p;\alpha ,\beta } & \le c_1 &
\int^1_{-1} (1-t)^{\alpha (1-p/2)-p/4} (1+t)^{\beta (1-p/2)-p/4}dt \\
 & = & c_2 \ = \ c_2 \parallel p_j^{(\alpha ,\beta )} \parallel_{2;\alpha
 ,\beta }< \infty
\end{eqnarray*}
and
\[
\parallel \sup_j |p_j^{(\alpha ,\beta )} | \parallel_{2; \alpha ,\beta }
\le c_3 (\int^1_{-1} (1-t^2)^{-1/2} dt)^{1/2} =c_4.
\]

Thus $p=2$ and $X$ is a Hilbert space by Proposition 4.
{\eproof}
The unconditionality of the Haar system in $L_p(0,1;X)$, if $1<p<\infty$
and $X$ is an UMD-space, shows that Proposition 4 does not hold without
conditions being imposed on the system $(p_j)$ as done in (i), (ii)
there.

For the proof of Propositions 5 and 7, we need the following result due
to Lindenstrauss and Pe\l czy\'nski \cite[proof of Theorem 4.2]{kn:lp} and
Olevskii \cite{kn:o}.

\begin{theorem}
Let $1 \le p \le \infty $ and $(p_n)_{n \in \Bbb N}$ be a basis of $L_p
(0,1)$. Let $(h_j)_{j \in \Bbb N}$ denote the Haar system on
$[0,1]$, normalized by $\parallel h_j \parallel_p =1$. For any $0 <
\delta < 1$ there is a block basis sequence $(z_j)_{j \in \Bbb N}$ of
$(p_n)_{n \in \Bbb N}$ such
that for every $N \in \Bbb N$ there is a measure preserving automorphism
$\varphi_N :[0,1] \rightarrow [0,1]$ with 
\[
\sum^N_{j=1} \parallel z_j^{\ast} \parallel_{p'} \parallel h_j o \varphi_N -z_j
\parallel_p \le \delta,
\]
where $(z_j^{\ast}) \subset L_{p'} (0,1)$ is biorthogonal to $(z_j) \subset
L_p (0,1)$.
\end{theorem}

Thus for some increasing sequence $(m_j)_{j \in \Bbb N_0}$ of integers
and scalars $(a_n)_{n \in \Bbb N}$
\[
z_j = \sum^{m_j}_{n=m_{j-1}+1} a_n p_n
\]
is in the above sense close to the Haar system.

{\bproof} {\bf of Proposition 5.} We will show that the Haar system is
unconditional in $L_p (0,1;X)$. Then, by Maurey \cite{kn:ma} $X$ has to
be a UMD-space, using also the results of \cite{kn:bo} and
\cite{kn:bu}.

Let $N \in \Bbb N$, $0 < \delta < 1$ and $x_1, \cdots , x_N \in X$. Let
$\varphi_N$ be as in the theorem and put $g_j :=h_j o \varphi_N$. Since
$(p_n)_{n \in \Bbb N}$ is unconditional in $L_p (0,1;X)$ by assumption,
so is the block basic sequence $(z_j)_{j \in \Bbb N}$. Hence for any
sequence of signs $(\varepsilon_j), \varepsilon_j \in \{ +1,-1 \}$, 

\begin{eqnarray*}
(\int^1_0 \parallel \sum^N_{j=1} \varepsilon_j h_j (t) x_j \parallel^p
dt)^{1/p} & = & (\int^1_0 \parallel \sum^N_{j=1} \varepsilon_j g_j (t) x_j
\parallel^p dt )^{1/p} \\
 & \le & (\int^1_0 \parallel \sum^N_{j=1} \varepsilon_j z_j (t) x_j
 \parallel^p dt)^{1/p} + (\int^1_0 \parallel \sum^N_{j=1} \varepsilon_j
 (g_j (t) - z_j (t))x_j \parallel^p dt )^{1/p} \\
 & \le & K (\int^1_0 \parallel \sum^N_{j=1} z_j (t) x_j \parallel^p
 dt)^{1/p} + \sum^N_{j=1} \parallel g_j - z_j \parallel_p \parallel x_j
 \parallel \\
 & \le & (K+\delta ) (\int^1_0 \parallel \sum^N_{j=1} z_j (t) x_j
 \parallel^p dt )^{1/p},
\end{eqnarray*}
using that $x_j = <z_j^{\ast} , \sum^N_{k=1} z_k x_k >$ and hence 
\[
\parallel x_j \parallel \le \parallel z_j^{\ast} \parallel_p (\int^1_0
\parallel \sum^N_{k=1} z_k (t) x_k \parallel^p dt)^{1/p}.
\]

The constant $K$ is independent of $N$, $(x_j)$ and $(\varepsilon_j)$. The
chain of inequalities can be reversed with all $\varepsilon_j = +1$ to
find

\begin{eqnarray*}
(\int^1_0 \parallel \sum^N_{j=1} \varepsilon_j h_j (t) x_j \parallel^p
dt)^{1/p} & \le & (K + \delta ) (\int^1_0 \parallel \sum^N_{j=1} z_j (t)
x_j \parallel^p dt )^{1/p} \\
 & \le & (K+\delta )(1+ \delta ) (\int^1_0 \parallel \sum^N_{j=1} h_j (t)
 x_j \parallel^p dt)^{1/p},
\end{eqnarray*}
i.e. the Haar system is unconditional.
{\eproof}
A similar procedure is used in the

{\bproof} {\bf of Proposition 7:} Let $N \in \Bbb N$, $x_1, \cdots x_N \in
X$ and $0 < \delta < 1$. With the same notation as in the previous
proof,
\[
(\int^1_0 \parallel \sum^N_{j=1} h_j (t) x_j \parallel^2 dt)^{1/2} \le
(1+ \delta ) (\int^1_0 \parallel \sum^N_{j=1} z_j (t) x_j \parallel^2
dt)^{1/2}.
\]
with $z_j = \sum^{m_j}_{n=m_{j-1}+1} a_n p_n$ and $(m_j)_{j \in \Bbb N_0}
\subset \Bbb N$ increasing. Note that 
\[
(\sum^{m_j}_{n=m_{j-1}+1} |a_n|^2)^{1/2} = \parallel z_j \parallel_2 \le
(1+\delta ) \parallel h_j \parallel_2 = (1+\delta ).
\]

Thus, using that $X$ has $(p_n)$-type 2, there is $K$ independent of
$x_1, \cdots x_N \in X$ such that

\begin{eqnarray*}
(\int^1_0 \parallel \sum^N_{j=1} h_j (t) x_j \parallel^2 dt)^{1/2} & \le &
(1+\delta ) (\int^1_0 \parallel \sum^N_{j=1} (\sum^{m_j}_{n=m_{j-1}+1}
a_n p_n (t))x_j \parallel^2 dt)^{1/2} \\
 & \le & K(1+\delta )(\sum^N_{j=1} \sum^{m_j}_{n=m_{j-1}+1} |a_n|^2 \parallel
 x_j \parallel^2 )^{1/2} \\
 & \le & K(1+\delta )^2 (\sum^N_{j=1} \parallel x_j \parallel^2 )^{1/2}.
\end{eqnarray*}

This shows that $X$ has ``Haar-type 2'' which directly implies type 2
since the Rademacher functions form a block basis of the Haar functions,
\[
r_k = \sum^{n_k}_{j=n_{k-1}+1} t_j h_j  \ , \ \sum^{n_k}_{j=n_{k-1}+1}
|t_j|^2 =1.
\]
where $(n_k)_{k \in \Bbb N_0} \subset \Bbb N$ is a suitable increasing 
sequence. Hence for any sequence $y_1, \cdots , y_{\ell} \in X$

\begin{eqnarray*}
(\int^1_0 \parallel \sum^{\ell}_{k=1} r_k (t) y_k \parallel^2
dt)^{1/2} & = & (\int^1_0 \parallel \sum^{\ell}_{k=1}
(\sum^{n_k}_{j=n_{k-1}+1} t_j h_j (t))y_k \parallel^2 dt)^{1/2} \\
 & \le & K(1+\delta )^2 (\sum^{\ell}_{k=1}(\sum^{n_k}_{j=n_{k-1}+1}
 |t_j|^2 ) \parallel y_k \parallel )^{1/2} \\
 & = & K(1+ \delta )^2 (\sum^{\ell}_{k=1} \parallel y_k \parallel^2
 )^{1/2}.
\end{eqnarray*}
{\eproof}
\vspace{1 ex}
We note that if conversely $X$ has type 2, one has as estimate of the
Rademacher against the Haar mean, i.e. there is a constant $C$ such that for
all $\ell \in \Bbb N$ and all \mbox{$(y_k)^{\ell}_{k=1} \subseteq X$:}
\[
\int^1_0 \parallel \sum^{\ell}_{k=1} r_k(t) y_k \parallel^2 dt \le
C \int^1_0 \parallel \sum^{\ell}_{k=1} h_k (t) y_k \parallel^2 dt.
\]

Indeed, this statement is equivalent to the existence of a constant $C_1$, so 
that for all $\ell \in \Bbb N$ and all $f \in L_2 (0,1;X)$
\[
(\int^1_0 \parallel \sum^{\ell}_{k=1} r_k(t)<f,h_k> \parallel^2)^{1/2} \le
C_1 (\parallel f \parallel_2)^{1/2}
\]

Let $f \in L_2(0,1;X)$ be of the form $f:= \sum^n_{j=1} x_j \otimes f_j$, 
where \mbox{$(f_j)^n_{j=1} \subseteq L_2 (0,1)$} is a finite sequence of
normalized, mutually disjointly supported functions and \mbox{$(x_j)^n_{j=1}
\subseteq X$.} Since the set of such functions is a dense subspace of
$L_2 (0,1;X)$ it suffices to prove the inequality for those.

Let $(\gamma_j)$ be a sequence of standard independent $N(0,1)$ Gaussian 
variables on a probability space $(\Gamma, \nu)$ and let M be the type 2
constant of $X$. Since for every $x^{\ast} \in X^{\ast}$
\begin{eqnarray*}
\sum^{\ell}_{k=1} |x^{\ast}(<f,h_k>)|^2 & = & \sum^{\ell}_{k=1} |<x^{\ast}f,
h_k>|^2 \\
 & \le & \int^1_0 |x^{\ast}(f(t))|^2 dt \\
 & = & \sum^n_{j=1} |x^{\ast}(x_j)|^2,
\end{eqnarray*}
it follows from the unitary invariance of the $\gamma_j$'s (see e.g. 
\cite{kn:mp}) that
\[
\int_{\Gamma} \parallel \sum^{\ell}_{k=1} \gamma_k (s) <f,h_k> \parallel ^2 
d{\nu}(s) \le \int_{\Gamma} \parallel \sum^n_{j=1} \gamma_j (s) x_j 
\parallel ^2
d{\nu}(s).
\]

Combining this with the fact that since $X$ is of type 2 the Rademacher and the
Gauss means are $K$-equivalent for a suitable $K$, we obtain
\begin{eqnarray*}
(\int^1_0 \parallel \sum^{\ell}_{k=1} r_k(t)<f,h_k> \parallel ^2 dt)^{1/2}
& \le & K (\int_{\Gamma} \parallel \sum^{\ell}_{k=1} \gamma_k (s)<f,h_k>
\parallel ^2 d\nu (s))^{1/2} \\
& \le & K (\int_{\Gamma} \parallel \sum^n_{j=1} \gamma_j (s)x_j \parallel ^2 
d\nu (s))^{1/2} \\
& \le & KM (\sum^n_{j=1} \parallel x_j \parallel ^2)^{1/2}  \ = \ \parallel f
\parallel _2,
\end{eqnarray*}
which proves the claim.

The partial converse of Proposition 7 follows easily:

{\bproof} {\bf of Proposition 8:} Let \mbox{$N \in \Bbb N$}, 
$x_1, \cdots x_N \in
X$. By the unconditionality assumption on the $(p_n)$-system in
$L_2(0,1;X)$ and the type property of $X$ there are constants $c_1, c_2$
independent of $N$ and $x_1, \cdots x_N \in X$ such that

\begin{eqnarray*}
(\int^1_0 \parallel \sum^N_{j=1} p_j (t) x_j \parallel^2 dt)^{1/2} & \le &
c_1 (\int^1_0 \int^1_0 \parallel \sum^N_{j=1} r_j (s) p_j (t) x_j
\parallel^2 ds dt)^{1/2} \\
 & \le & c_2 (\int^1_0 \sum^N_{j=1} |p_j (t) |^2 \parallel x_j \parallel^2
 dt )^{1/2} \\
 & = c_2 & (\sum^N_{j=1} \parallel x_j \parallel^2 )^{1/2}.
\end{eqnarray*}
Hence $X$ has $(p_n)$-type 2.
{\eproof}
\vspace{1 ex}
We still have to show that in certain cases $X$ is isomorphic to a
Hilbert space provided that it only has $(p_n)$-type 2. This will be
another application of the interpolation inequalities of Theorem 1.

{\bproof} {\bf of Proposition 6:} For $n \in \Bbb N$, let
$(t_j)^{n+1}_{j=1}$ denote the Gaussian quadrature weights. The
$(n+1)\times (n+1)$ matrix $A_n = (a_{ij})$ defined by
$a_{jk}:=\sqrt{\lambda_j} \, p_k^{(\alpha ,\beta )} (t_j), j=1, \cdots ,
n+1, k=0, \cdots ,n$ is orthogonal since by Gaussian quadrature for
$k,{\ell} \in {0, \cdots ,n}$

\begin{eqnarray*}
\delta_{k{\ell}} & = & \int^1_{-1} p_k^{(\alpha ,\beta )} (t) p_{\ell
}^{(\alpha ,\beta )} (t)w_{\alpha \beta}(t) dt \\
 & = & \sum^{n+1}_{j=1} \lambda_j p_k^{(\alpha , \beta)}(t_j) p_{\ell}
^{(\alpha, \beta)}(t_j) \\
 & = & \sum^{n+1}_{j=1} a_{jk} a_{j{\ell}}.
\end{eqnarray*}

Since the measure space $(I,w_{\alpha \beta })$ is equivalent to
$(0,1)$, we know from Proposition 7 that $X$ has type 2. We will now
show that $X$ also has cotype 2 and hence by Kwapie\'n \cite{kn:kw} is
isomorphic to a Hilbert space. We use Theorem 1 (a) to discretize the
notion of $(p_n^{(\alpha ,\beta )})$-type 2 and reverse the inequality
using the orthogonality of the matrix $A_n$ appearing in this way: By
Theorem 1 and the $(p_n^{(\alpha ,\beta )})$-type 2 property there are
$c_1, c_2$ such that for any $n \in \Bbb N$ and $x_0, \cdots , x_n \in
X$
\begin{eqnarray}
(\sum_{j=1}^{n+1} \parallel \sum_{k=0}^n a_{jk} x_k \parallel^2)^{1/2} &
= & (\sum_{j=1}^{n+1} \lambda_j \parallel \sum^n_{k=0} p_k^{(\alpha ,\beta
)} (t_j) x_k \parallel^2 )^{1/2} \nonumber \\
 & \le & c_1 (\int^1_{-1} \parallel \sum_{k=0}^n p_k^{(\alpha ,\beta )}
 (t) x_k \parallel^2 w_{\alpha \beta } (t) dt)^{1/2} \nonumber \\
 & \le & (\sum^n_{k=0} \parallel x_k \parallel^2 )^{1/2}. 
\end{eqnarray}

Since $A^{-1}_n = {A_n}^t$, we can invert (27) easily: starting with
arbitrary $y_1, \cdots , y_{n+1} \in X$ and applying (27) to $x_k:=
\sum^{n+1}_{{\ell} =1} a_{{\ell} k} y_{\ell}$, $k \in {0, \cdots ,n}$, we
find
\[
(\sum^{n+1}_{j=1} \parallel y_j \parallel^2 )^{1/2} \le c_2
(\sum^n_{k=0} \parallel \sum_{{\ell}=1}^{n+1} a_{{\ell}k} y_{\ell}
\parallel^2 )^{1/2}.
\]

If the $A_n$ were symmetric, this and Theorem 1 (a) would yield the cotype 2
property. However, $A_n \ne {A_n}^t$, in general. To prove the cotype
2 property, we replace $y_j$ by $r_j(s)y_j$ and apply the contraction
principle to find

\begin{eqnarray*}
(\sum_{j=1}^{n+1} \parallel y_j \parallel^2 )^{1/2} & \le & c_2
(\sum^n_{k=0} \int^1_0 \parallel \sum_{{\ell}=1}^{n+1} a_{{\ell}k}
r_{\ell}(s)y_{\ell} \parallel^2 ds )^{1/2} \\
 & \le & c_3 (\sum^n_{k=0} \sup_{{\ell} \le n+1} |a_{{\ell}k} |^2
 )^{1/2} (\int^1_0 \parallel \sum^{n+1}_{{\ell}=1} r_{\ell} (s) y_{\ell}
 \parallel^2 ds)^{1/2}.
\end{eqnarray*}

Hence $X$ will have cotype 2 provided that $\sum^n_{k=0} \sup_{{\ell}
\le n+1} |a_{{\ell}k} |^2$ is uniformly bounded in $n \in \Bbb N$.
This is correct since by (4) and (10) e.g. if ${\ell} \le n/2$
\[
\sqrt{\lambda_{\ell}} \sim {\ell}^{\alpha + 1/2}/n^{\alpha+1} \ , \  
|p_k^{(\alpha ,\beta )} (t_{\ell}) | \stackrel{<}{\sim}(n/{\ell})^{\alpha +1/2} 
\]
\[
|a_{{\ell}k}| = \sqrt{\lambda_{\ell}} \, |p_k^{(\alpha , \beta)}(t)| 
\stackrel{<}{\sim}
n^{-1/2}
\]   

The case of ${\ell} > n/2$ is similar. By the result of Kwapie\'n, used
earlier,  $X$ is isomorphic to a Hilbert space.
{\eproof}

\begin{tabbing}

xxxxxxxxxxxxxxxxxxxxxxxxxxxxxxxx \= xxxxxxxxxxxxxxxxxxxxxxxxxxxxxxxxxxx
\kill

Hermann K\"onig \> Niels J\o{}rgen Nielsen \\
Mathematisches Seminar \> Matematisk Institut \\
Universit\"at Kiel \> Odense Universitet \\
Ludewig Meyn Str. 4 \> Campusvej 55 \\
23 Kiel \> 5230 Odense M \\
Germany \> Denmark 
\end{tabbing}
\end{document}